\documentclass[a4paper]{article}
\usepackage[utf8]{inputenc}
\usepackage[T2A]{fontenc}

\usepackage{amssymb}
\usepackage{amsmath}
\usepackage[all]{xy}
\usepackage[CJKbookmarks=true]{hyperref}

\title{Cauchy Formulas and Billey's Formulas for Generalized Grothendieck polynomials}
\author{XIONG Rui}
\date{}

\newtheorem{Th}{Theorem}[section]
\newtheorem{Lemma}[Th]{Lemma}
\newtheorem{Coro}[Th]{Corollary}
\newtheorem{Prop}[Th]{Proposition}

\newcounter{term}[section]
\def\theterm{\thesection.\arabic{term}}
\def\term{\futurelet\next\termfoo}
\def\termfoo{\ifx\next[%
    \expandafter\termAux%
\else%
    \expandafter\termaux%
\fi%
}
\def\termAux[#1]{\refstepcounter{term}%
\paragraph{\theterm. #1}}
\def\termaux{\refstepcounter{term}%
\paragraph{\theterm}}

\begin{document}

\def\CH{\operatorname{CH}}
\def\codim{\operatorname{codim}}
\def\pt{\operatorname{\mathsf{pt}}}

\def\ru{\hline\rule[-.5pc]{0pc}{1.5pc}}
\def\Ad{\operatorname{Ad}}

\def\GL{\operatorname{GL}}
\def\Fl{\operatorname{\mathcal{F}\ell}}

\def\NH{\operatorname{NH}}

\def\ppi{\hat{\pi}}

\def\longdash{\mathord-\!\!\!\mathord-\!\!\!\mathord-\!\!\!\mathord-}
\def\longDash{\mathord=\!\!\!\mathord\Rightarrow\!\!\!\mathord=\!\!\!\mathord=}
\def\longDdash{\mathord\equiv\!\!\!\mathord\Rrightarrow\!\!\!\mathord\equiv\!\!\!\mathord\equiv}

\setlength{\unitlength}{1pc}
\def\labe#1{%
\begin{array}{@{}l@{}}\vspace{-0.3pc}\begin{picture}(1,1)
  \put(0.25,0.2){$#1$}
\end{picture}\end{array}}
\def\empt{%
\begin{array}{@{}l@{}}\vspace{-0.3pc}\begin{picture}(1,1)
  \put(0,0){\line(1,0){1}}\put(0,0){\line(0,1){1}}
  \put(0,1){\line(1,0){1}}\put(1,0){\line(0,1){1}}
\end{picture}\end{array}}
\def\pipe{%
\begin{array}{@{}l@{}}\vspace{-0.3pc}\begin{picture}(1,1)
  \put(0,0){\line(1,0){1}}\put(0,0){\line(0,1){1}}
  \put(0,1){\line(1,0){1}}\put(1,0){\line(0,1){1}}
  \linethickness{1pt}
  \qbezier(0,0.5)(0.5,0.5)(1,0.5)\qbezier(0.5,0)(0.5,0.5)(0.5,1)
\end{picture}\end{array}}
\def\bump{%
\begin{array}{@{}l@{}}\vspace{-0.3pc}\begin{picture}(1,1)
  \put(0,0){\line(1,0){1}}\put(0,0){\line(0,1){1}}
  \put(0,1){\line(1,0){1}}\put(1,0){\line(0,1){1}}
  \linethickness{1pt}
  \qbezier(0,0.5)(0.5,0.5)(0.5,1)\qbezier(0.5,0)(0.5,0.5)(1,0.5)
\end{picture}\end{array}}
\def\abou{%
\begin{array}{@{}l@{}}\vspace{-0.3pc}\begin{picture}(1,1)
  \qbezier(0,0.5)(0.5,0.5)(0.5,1)
\end{picture}\end{array}}

\def\Bpipe{%
\begin{array}{@{}l@{}}\vspace{-0.3pc}\begin{picture}(1,1)
  \put(0,0){\line(1,0){1}}\put(0,0){\line(0,1){1}}
  \put(0,1){\line(1,0){1}}\put(1,0){\line(0,1){1}}
  \linethickness{1pt}
  \qbezier(0,0.5)(0.5,0.5)(0.5,1)\qbezier(0.5,0)(0.5,0.5)(1,0.5)
  \qbezier(0.2,0.2)(0.5,0.5)(0.8,0.8)
\end{picture}\end{array}}

\def\Bbump{%
\begin{array}{@{}l@{}}\vspace{-0.3pc}\begin{picture}(1,1)
  \linethickness{1pt}
  \qbezier(0,0.5)(0.5,0.5)(0.5,1)\qbezier(0.5,0)(0.5,0.5)(1,0.5)
\end{picture}\end{array}}
\def\abou{%
\begin{array}{@{}l@{}}\vspace{-0.3pc}\begin{picture}(1,1)
  \qbezier(0,0.5)(0.5,0.5)(0.5,1)
\end{picture}\end{array}}

\def\wt{\operatorname{\mathtt{wt}}}

\maketitle

\noindent\textbf{Abstract}\quad 
We study the generalized double $\beta$-Grothendieck polynomials for all types.  
We study the Cauchy formulas for them. 
Using this, we deduce the K-theoretic version of the comodule structure map $\alpha^*: K(G/B)\to K(G)\otimes K(G/B)$ induced by the group action map for reductive group $G$ and its flag variety $G/B$. 
Furthermore, we give a combinatorial formula to compute the localization of Schubert classes as a generalization of Billey's formula. 

\bigbreak
I would politely express my gratitude to 
Victor Petrov, 
Neil JiuYu Fan,  
Peter Long Guo and 
for discussion.



\tableofcontents

\section{Main Results}

\term
Let $G$ be a connected complex reductive group, 
and $B$ its Borel subgroup.
The homogenous variety $G/B$ is called the \emph{flag variety} of $G$.

Let $W=N_G(T)/T$ be the Weyl group and $\ell$ the standard length function. 
For elements $u,v,w\in W$, we write \emph{the reduced decomposition} $w=u\odot v$ if $w=uv$ and $\ell(w)=\ell(u)+\ell(v)$. 

We also introduce the product $*$ over $W$ which is uniquely characterized by 
$$\forall i\in\mathbb{I},\,\,s_i * s_i=s_i\qquad uv=u\odot v\Longrightarrow uv=u*v,$$
with $\{s_i:i\in \mathbb{I}\}$ the set of simple reflections. 
Geometrically, $\Sigma_u\cdot \Sigma_{v}=\Sigma_{u*v}$ for any $u,v\in W$ by the Tits system. 

\term 
For an element $w\in W$, we define the \emph{lower Schubert variety} $\Sigma_w$ to be the Zariski closure of $BwB/B\subseteq G/B$, 
and the \emph{(upper) Schubert variety} $\Sigma^w$ to be the Zariski closure of  $B^-wB/B\subseteq G/B$, where $B^-$ is the opposite Borel subgroup. 
It is well-known \cite{SpringerLinear} that the
$\dim_{\mathbb{C}} \Sigma_w=\codim_{\mathbb{C}} \Sigma^w=\ell(w)$. 
Furthermore, $\Sigma_w$ and $\Sigma^w$ are both isomorphic to affine linear spaces.

\term Let $\beta$ be a parameter. Denote $K(-)=K^\beta(-)$ the $\beta$-Grothendieck group. 
It is the oriented cohomology theory universal with respect to the multiplicative form group law
$$x\oplus_\beta y= x+y-\beta xy. $$
See for example \cite{levine2007algebraic}. After specialization at $\beta=1$, we get the usual Grothendieck group $K_0$. 

We also use the equivariant form $K_B(-)$. For the case $\beta=1$, the basic definition can be found in
\cite{chriss2009representation}. 
The general case can be easily established similarly. 

Over the flag variety $G/B$, we denote $\mathcal{O}^w$ the push forward of $\mathcal{O}_{\Sigma^w}$ to $G/B$. We denote $[\mathcal{O}^w]$ the class of it in $K(G/B)$ and $[\mathcal{O}^w]_{B}$ the equivariant analogy. 
Then by the affine stratification theorem, 
$K(G/B)$ is freely generated by $[\mathcal{O}^w]$ over $\mathbb{Z}[\beta,\beta^{-1}]$, and $K_B(G/B)$ is freely generated by $[\mathcal{O}^w]_B$ over $K_B(\pt)$. 

\term Let $K^{\text{top}}(-)$ be the complex topological K-theory. We have a map 
$\eta_X:K^{\beta=1}(X)\to K^{\text{top}}(X)$. 
By topological K-theory, $\eta_{G/B}$ is an isomorphism. 
The equivariant version is also true but  after completion, 
by the Atiyah completion theorem \cite{atiyah1969equivariant} and
the Atiyah--Hirzebruch spectral sequence.

\term We have the following analogy of the main result of \cite{xiong2020comodule}. 

\begin{Th}\label{Descriptionofalpha}The map induced by the left action $\alpha:G\times G/B\to G/B$
$$\alpha^*: K(G/B)\longrightarrow K(G)\otimes K(G/B)$$
is given by
$$[\mathcal{O}^w] \longmapsto
\sum\nolimits_{w=u*v} (-\beta)^{\ell(v)+\ell(u)-\ell(w)} \pi^*[\mathcal{O}^{u}]\otimes [\mathcal{O}^{v}].$$
The topological K-theory has the same formula after specialization at $\beta=1$. 
\end{Th}

\term Other than consider the Borel construction, we think $G$ with $B$ acting on both sides, and 
$G\times G$ with $B$ acting on left, right, and middle
$$(b_1,b_2,b_3)\cdot (g,h)=(b_1gb_2^{-1} ,b_2hb_3^{-1}). $$
Then the multiplication map $\mu:G\times G\to G\times B$ is ($B\times B\times B$)-equivariant. 

We denote the obvious map 
$$K_{B\times B\times B}(G\times G)
\begin{cases}
\stackrel{\pi_1}{\longleftarrow}K_{G\times B\times B}(G\times G)\cong K_{B\times B}(G)\\
\stackrel{\pi_2}{\longleftarrow}K_{B\times B\times G}(G\times G)\cong K_{B\times B}(G)\\
\end{cases}$$

\term We have the following analogy of the equivariant version of \cite{xiong2020comodule}. 

\begin{Th}\label{Descriptionofmu}
The map induced by the multiplication $\mu$ 
$$\mu^*: K_{B\times B}(G;\mathbb{Q})\longrightarrow K_{B\times B\times B}(G\times G;\mathbb{Q})$$
is given by 
$${}[\mathcal{O}^w]_B \longmapsto
\sum\nolimits_{w=u*v} (-\beta)^{\ell(v)+\ell(u)-\ell(w)} \pi_1^*[\mathcal{O}^{u}]_B\cdot \pi_2^*[\mathcal{O}^{v}]_B.$$
The topological K-theory has the same formula after specialization at $\beta=1$. 
\end{Th}

\term 
Denote 
$R_T^{\beta}=\mathbb{Q}[e^{\beta\Lambda}]\otimes \mathbb{Q}[\beta,\beta^{-1}]$, 
obtained by adjoint a new variable $\beta$ to the group ring of $\beta\Lambda\cong \Lambda$. 
Let $R_G=(R_T)^W$ be the invariant subalgebra of $R_T$. 
Then Borel's presentation still holds  
$$K(G/B;\mathbb{Q})=R_T\otimes_{R_G}\mathbb{Q},\qquad 
K_B(G/B;\mathbb{Q})=R_T\otimes_{R_G}R_T. $$
Actually, the above map is compatible with the Borel's presentation of cohomology under the Chern character.

To be precise, $\mathcal{O}(\lambda)$ is presented by $e^{\beta\lambda}$. Note that its first Chern class is $\frac{e^{\beta\lambda}-1}{\beta}$ by definition. 

\term Then the Schubert class $[\mathcal{O}^w]\in K(G/B;\mathbb{Q})$ corresponds to some element in $R_T\otimes_{R_G}\mathbb{Q}$, which we will denoted by $\mathfrak{G}_w(X)$ and call it the \emph{generalized $\beta$-Grothendieck polynomial}. 
On the other hand, the equivariant analogy for $[\mathcal{O}^w]_B\in H_B^\bullet(G/B;\mathbb{Q})$ is denoted by $\mathfrak{G}_w(X,T)$ and is called \emph{the generalized double $\beta$-Grothendieck polynomial}. Then $\mathfrak{G}(X,1)=\mathfrak{G}(X)$. To be precise, $T\mapsto 1$ is the map sending all $1\otimes e^{\beta \lambda}$ to $1$.

\term We have the following analogy of Cauchy formulas proved in \cite{xiong2020comodule}. 

\begin{Th}\label{CauchyformulaforGrothendieck}We have 
$$\mathfrak{G}_w(X,Z)=\sum\nolimits_{u*v=w}(-\beta)^{\ell(u)+\ell(v)-\ell(w)} \mathfrak{G}_v(X,Y)\mathfrak{G}_u(Y,Z)$$
in $R_T\otimes_{R_G}R_T\otimes_{R_G}R_T$. 
\end{Th}

\term In the case $G=\GL_n$, the generalized Grothendieck polynomials has a stable choice as in the cohomology case. 
This coincides with the definition of \cite{fomin1993yang} after replacing $\beta$ by $-\beta$. 
In particular, when $\beta=1$, it recovers the usual Grothendieck polynomials. 
We will shortly review the combinatorics in the appendix.

\term Actually we introduce a generalized double dual $\beta$-Grothendieck polynomial with plenty of good properties. They have an inversion formula (Theorem \ref{inversionFo1} and Theorem \ref{inversionFo2}). The Cauchy formula (Theorem \ref{StarCauchyFor} and Theorem \ref{StarCauchyFor2}). They are the dual basis of the Demazure operators (Theorem \ref{DualFor1} and Theorem \ref{DualFor2}).

\term Let $T$ be the maximal torus $T$ of $G$ contained in $B$. 
For the point $w\cdot B/B$, there is a localization map 
$$K^\beta_B(G/B)\cong K^\beta_T(G/B)\longrightarrow K^\beta_T(wB/B)\cong K_T^\beta(\pt). $$
We give a combinatorial description of the image of  $[\mathcal{O}^u]$ under this localization map.  

\term Assume that $w=s_{i_1}s_{i_2}\cdots s_{i_r}$ (not necessarily reduced). 
For any $j=1,\ldots,r$, denote 
$$d_{j}=s_{i_1}\cdots s_{i_{j-1}}\bar{a}_j, $$
where $a_i=\frac{1-e^{-\beta\alpha_i}}{\beta}$, and $\bar{a}_i=\ominus_\beta a_i:=\frac{-a_i}{1-\beta a_i}$. 
For any subset $J=\{j(1)<\cdots<j(k)\}\subseteq \{1,\ldots,r\}$, denote 
$$w(J)=s_{i_{j(1)}}*\cdots*s_{i_{j(k)}}. $$

\begin{Th}\label{GeoBilleyformula}
The image of $[\mathcal{O}^u]_B$ under this localization map to $wB/B$ is given by $$\mathfrak{G}_u(wx,x)=\sum_{w(J)=u}(-\beta)^{|J|-\ell(u)}\prod_{j\in J} d_j\in K_T^\beta(\pt). $$
Note that this expression does not depends on the choice of the decomposition of $w$. 
\end{Th}

\term Taking $\beta \to 0$, this recovers \emph{Billey's formula} \cite{billey1999kostant}. 
But actually, this formula has been discovered in \cite{andersen1994representations}. 
See \cite{tymoczko2013billeys} for remarks and applications. 
This formula also generalizes the Buch--Rim\'anyi formula \cite{buch2004specializations} in $A$-Type Grothendieck polynomial ($\beta=1$).

\section{Demazure Operators and Localization}

To simplify notations, all Grothendieck groups are of coefficients in $\mathbb{Q}$ in this section. 

\term 
For the standard parabolic subgroup $P_i=B\cup Bs_iB$, the group $K_{P_i}^\bullet(\pt)$ can be computed to be $R_T^{s_i}$ the subalgebra of $R_T$ fixed by $\{1,s_i\}$. 
Note that the $G$-equivariant K-theory $K_G(G/H)=K_H(\pt)$ for any closed subgroup $H$. 
We introduce the Demazure operator $\pi_i:R_T\to R_T$ to be the composition of
$$\begin{array}{c}
\sigma_*:R_T=K_G(G/B)\longrightarrow H_G^{\bullet-2}(G/P_i)=R_T^{s_i},\\ [1ex]
\sigma^*:R_T^{s_i}=K_G(G/P_i)
\longrightarrow K_B(G/B)=R_T
\end{array}$$
where $\sigma:G/B\to G/P_i$ is the natural projection, inducing the Gysin push forward $\sigma_*$ and pull back $\sigma^*$.

\term By definition, the $B$-equivariant cohomology Demazure operator $R_T\otimes_{R_G}R_T\to R_T\otimes_{R_G}R_T$, the composition
$$\begin{array}{c}
\sigma_*:R_T\otimes_{R_G} R_T=K_B(G/B)\longrightarrow H_B^{\bullet-2}(G/P_i)=R_T^{s_i}\otimes_{R_G} R_T,\\ [1ex]
\sigma^*:R_T^{s_i}\otimes_{R_G} R_T=K_B(G/P_i)
\longrightarrow K_B(G/B)=R_T\otimes_{R_G} R_T
\end{array}$$
is given by $\pi_i\otimes 1$.

\term It is well-known that over $K_B(G/B)$,
$$\ell(ws_i)=\ell(w)-1\quad \Longrightarrow\quad \pi_i[\mathcal{O}^w]_B=[\mathcal{O}^{ws_i}]_B.$$
In terms of generalized double $\beta$-Grothendieck polynomials, over $R_T^\beta\otimes_{R_G^\beta}R_T^\beta$, 
$$\ell(ws_i)=\ell(w)-1\quad \Longrightarrow\quad \pi_i^X\mathfrak{G}_w(X,Y)=\mathfrak{G}_{ws_i}(X,Y)$$
where $\pi_i^X=\pi_i\otimes 1$. 
At the present stage, we cannot obtain the case $\ell(ws_i)=\ell(w)+1$ directly, but it will be described later. 

\term 
To obtain the K-theory version of Demazure operators, we use the Grothendieck--Riemann--Roch theorem for $\beta$-Grothendieck group \cite{levine2007algebraic}. 
To be precise, the relative tangent bundle of $G/B\to G/P_i$ is $\mathcal{O}(\alpha_i)$ whose $\beta$-Todd class is $\frac{\beta\alpha_i}{1-e^{-\beta\alpha_i}}$. Thus we have 
$$\pi_i f=\frac{\frac{\beta\alpha_i}{1-e^{-\beta\alpha_i}}f-
s_i\big(\frac{\beta\alpha_i}{1-e^{-\beta\alpha_i}}f\big)}{\alpha_i}
=\beta \frac{\frac{\alpha_i}{1-e^{-\beta\alpha_i}}f-
\frac{-\alpha_i}{1-e^{\beta\alpha_i}}
s_if}{\alpha_i}=\beta\frac{f-e^{-\beta\alpha_i}s_if}{1-e^{-\beta\alpha_i}}. $$
Actually, when $\beta=1$, this is what Demazure originally obtained in \cite{demazure1974desingularisation}. 

Taking into account of the action of $\pi_i$ on $[\mathcal{O}^w]_B$, we can easily see that $\pi_i$ for $i\in \mathbb{I}$ satisfies the Braid relation. 
By direct computation, we see that $\pi_i^2=\beta\pi_i$. 
This implies that for each $w\in W$, we can define $\pi_w$ by $\pi_w=\pi_{i_1}\cdots \pi_{i_r}$ with $w=s_{i_1}\odot \cdots \odot s_{i_r}$ a reduce word decomposition.

\term As a corollary, 
$$\pi_i[\mathcal{O}^w]_B=\begin{cases}
{}[\mathcal{O}^{ws_i}]_B,& \ell(ws_i)=\ell(w)-1,\\
{}\beta[\mathcal{O}^w]_B,& \ell(ws_i)=\ell(w)+1.
\end{cases}$$
In terms of generalized double $\beta$-Grothendieck polynomials, $R_T^\beta\otimes_{R_G^\beta}R_T^\beta$, 
$$\pi_i^X\mathfrak{G}_w(X,Y)=\begin{cases}
{}\mathfrak{G}_{ws_i}(X,Y),& \ell(ws_i)=\ell(w)-1,\\
{}\beta\mathfrak{G}_{ws_i}(X,Y),& \ell(ws_i)=\ell(w)+1.
\end{cases}$$
Here $\pi_i^x=\pi_i\otimes 1$. 

\term Note that the lifting of $W$ to $G/B$ is exactly $(G/B)^T$, the $T$-fixed point of $G/B$. By the K-theoretic localization theorem \cite{chriss2009representation}, the map
$$K_B(G/B)\longrightarrow \bigoplus_{w\in W} K_B(w\cdot B/B)$$
induced by $(G/B)^T\subseteq G/B$ is injective. 
By a simple computation, the corresponding map 
$R_T\otimes_{R_G}R_T\longrightarrow 
\bigoplus\nolimits_{w\in W} R_T$
is given by $x\otimes y\mapsto ((wx)\cdot y)_{w\in W}$. 

\term For a closed $B$-subvariety $Y$ in $G/B$. 
Denote the image of $[\mathcal{O}_Y]_B$ under the above localization map to be $(a_w)_{w\in W}$. 
If the fixed point $w$ is not contained in $Y$, then $a_w=0$. 

For Schubert varieties, the fixed point $u\in \Sigma^w$ if and only $w_0u\leq w_0w$ in the Bruhat order \cite{SpringerLinear} where $w_0$ is the unique longest element of $W$. Equivalently, $w\leq u$. 
In terms of generalized double $\beta$-Grothendieck polynomials, 
$$\mathfrak{G}_{w}(uT,T)\neq 0\quad \Longrightarrow\quad 
w\leq u.$$
We will use the case
$\mathfrak{G}_{w}(T,T)=\begin{cases}
1, & w= \operatorname{id},\\
0, & \text{otherwise}.
\end{cases}$


\section{Leibniz Rules}

\term Let us denote the $\beta$-affine nil-Hecke algera $\NH^\beta(W)$ the algebra generated by left multiplications of elements of $R_T^\beta$ and $\pi_i$ with $i\in \mathbb{I}$ over $R_T^\beta$. 

We introduce the \emph{inverse $\beta$-Demazure operator} $\ppi_i=\pi_i-\beta\in \NH^\beta(W)$, that is
$$\ppi_if=
\beta  \frac{f-s_if}{e^{\beta \alpha_i}-1},$$
for $f\in R_T$. Consider the involution $D: f\mapsto e^{-\beta\rho} \overline{f}$, where $\overline{e^\lambda}=e^{-\lambda}$, and $\rho$ is the half sum of positive roots. 
By a direct computation 
$$\begin{array}{rl}
\pi_i(Df) & = \beta\dfrac{e^{-\beta \rho}\overline{f}-e^{-\beta\alpha_i}s_i(e^{-\beta \rho}\overline{f})}{1-e^{-\beta\alpha_i}}\\
& = \beta\dfrac{e^{-\beta \rho}\overline{f}-e^{-\beta \rho}s_i\overline{f}}{1-e^{-\beta\alpha_i}}
=e^{-\beta\rho} \beta \overline{\,\frac{f-s_if}{1-e^{\beta \alpha_i}}\,}
=-D(\ppi_if)
\end{array}$$
As a result, $\ppi_i =-D\pi_iD$, thus must satisfy braid relation with $\ppi_i^2=-\beta \ppi_i$. 

\term 
Finally, the computation
$$\begin{array}{rl}
\pi_i(fg)& = 
\beta\dfrac{f-e^{-\beta\alpha_i}s_if}{1-e^{-\beta\alpha_i}}g
+\beta e^{-\beta\alpha_i}\dfrac{g-s_ig}{1-e^{-\beta\alpha_i}}s_i f = (\pi_if) g+(s_if)(\ppi_i g)\\[2ex]
\ppi_i(fg)& = 
\beta\dfrac{f-s_if}{e^{\beta\alpha_i}-1}g
+\beta \dfrac{g-s_ig}{e^{\beta\alpha_i}-1}s_i f = (\ppi_if) g+(s_if)(\ppi_i g)\\
\end{array}$$
proves the Leibnize rule
$\begin{cases}
\pi_i(fg)=(\pi_if) g+(s_if)\ppi_i g,\\
\ppi_i(fg)=(\ppi_if) g+(s_if)\ppi_i g.
\end{cases}$.

\term The next lemma is a generalization of author's previous work \cite{xiong2020comodule}. 
Note that $\pi_i\ppi_i=\ppi_i\pi_i=0$, thus 
$$\forall w\in W, \qquad \ppi_i\pi_w=\begin{cases}
0 & \ell(s_iw)=\ell(w)-1\\
\pi_{s_iw}-\beta \pi_w & \ell(s_iw)=\ell(w)+1. 
\end{cases}$$

\begin{Lemma}[Top Leibniz rule]\label{TopLeibnize}We have the following 
$$\ppi_{w_0}((w_0f)g)=\sum\nolimits_{w\in W} (-1)^{\ell(w)}(\pi_{ww_0} f)\cdot  \pi_w g,$$
for all $f,g\in R_T$. 
\end{Lemma}
\textsc{Proof.} We can assume the left-hand-side to be 
$\sum_{u\in W}  (-1)^{\ell(u)} c_u\pi_{u} g$ with $c_u\in R_T$. 
We take $g$ to be a generic element in $R_T$. 
Apply $\ppi_i$, we get 
$$\begin{array}{l}
\displaystyle
\quad\ppi_i\bigg(\sum_{u\in W}  (-1)^{\ell(u)} c_u\pi_{u} g\bigg)\\
\displaystyle
=\sum_{u\in W}  (-1)^{\ell(u)} (\ppi_i c_u)\pi_{u} g
+\sum_{u\in W}  (-1)^{\ell(u)} \big(s_ic_u\big) \ppi_i\pi_{u} g\\
\displaystyle
=\sum_{u\in W} (-1)^{\ell(u)} (\ppi_ic_u)\cdot \pi_{u} g
+\sum_{\ell(s_iu)=\ell(u)+1} (-1)^{\ell(u)} (s_ic_u) (\pi_{s_iu}g-\beta \pi_ug) \\
\end{array}$$
As a result, 
$$-\beta c_w=
\begin{cases}
\ppi_i c_w -\beta s_ic_w & \ell(s_iw)=\ell(w)+1\\
\ppi_i c_w - s_ic_{s_iw}& \ell(s_iw)=\ell(w)-1\\
\end{cases}$$
As a result, 
$$\pi_i c_w=\begin{cases}
\beta s_ic_w, & \ell(s_iw)=\ell(w)+1,\\
s_ic_{s_iw},& \ell(s_iw)=\ell(w)-1.\\
\end{cases}$$
But $\pi_i c_w$ is already symmetric, thus
$$\pi_i c_w=\begin{cases}
\beta c_w, & \ell(s_iw)=\ell(w)+1,\\
c_{s_iw},& \ell(s_iw)=\ell(w)-1.\\
\end{cases}$$
By induction, we have 
$c_w=\pi_{ww_0} c_{w_0}$. 
Since 
$$\ppi_i(fg)=(s_if)\pi_ig+\big(\pi_if-\beta fg-\beta (s_if)\big)g, $$
it is easy to see $c_{w_0}=f$. 
\qquad Q.E.D. \bigbreak

\begin{Coro}\label{SymmetricTheorem}For any $f,g\in R_T$, 
$$\sum_{u\in W} (-1)^{\ell(u)}e^{\beta \rho} (\pi_{uw_0} f)\cdot  \pi_u g,$$
is symmetric under $W$.  
\end{Coro}
\textsc{Proof.} 
For $w\in W$, $s_i\rho=\rho- \alpha_i$, thus $\pi_ie^{-\beta\rho}=0$. 
So 
$$\forall f\in R_T, \qquad s_if=f\iff \pi_i(e^{-\beta\rho}f)=0$$
Then $\ppi_ig\cdot e^{\beta\rho}$ is always symmetric for any $g\in R_T$.
\qquad Q.E.D. \bigbreak

\term 
Let us denote the \emph{Demazure operator}
$\Pi_i\in \NH^\beta(W)$ by
$$\Pi_if:=\overline{\ppi_i\bar{f}}=-\beta\frac{f-s_if}{1-e^{-\beta\alpha_i}}. $$
They also satisfy the braid relations and $\Pi_i^2=-\beta \Pi_i$, thus we can write $\Pi_{w}$ for $w\in W$.
Note that 
$$e^{\beta\rho}\pi_w f
=\overline{e^{-\beta\rho}\overline{\pi_w f}}
=\overline{D(\pi_w f)}
=(-1)^{\ell(w)}\overline{\ppi_wDf}
=(-1)^{\ell(w)}\Pi_w \overline{Df}.$$

\begin{Coro}\label{SymmetricTheorem2}For any $f,g\in R_T$, 
$$\sum_{w\in W} (\Pi_{ww_0} f)\cdot  \pi_w g=\sum_{w\in W} (\Pi_{w} f)\cdot  \pi_{ww_0} g,$$
is symmetric under $W$.  
\end{Coro}

\term We can take the limit $\beta\to 0$ to get
$$\lim_{\beta\to 0} \pi_i=\partial_i,\qquad \lim_{\beta\to 0}\ppi_i=\partial_i,\qquad \lim_{\beta\to 0}\Pi_i=-\partial_i. $$
Here $\partial_i$ is the cohomology Demazure operator. 

\section{Cauchy Formulas}

\term The K-theory analogy is not absolutely direct from the proof of cohomology as in \cite{xiong2020comodule}, since $\pi_i\neq -\Pi_i$. 

\begin{Th}[Inversion Formula]\label{inversionFo1}
For any $f\in R_T^\beta$, we have the following identity in $R_T^\beta\otimes_{R_G^\beta}R_T^\beta$, 
$$f(Y)=\sum_{w\in W} \mathfrak{G}_w(X,Y)\Pi_{w^{-1}}^X f(X)$$
\end{Th}
\textsc{Proof.} Note that the right-hand-side is 
$$\pi_{w^{-1}w_0}^X\mathfrak{G}_{w_0}(X,Y)\Pi_{w^{-1}}^Xf(X)$$
which is symmetric in $X$ by Corollary \ref{SymmetricTheorem2}, so that we can take $X=Y$. 
\qquad Q.E.D. \bigbreak

\begin{Coro}[Dual Basis]\label{DualFor1}The operator $\Pi_{w^{-1}}$ and $\mathfrak{G}_w$ are dual to each other, that is, 
$$\Pi_{u^{-1}}^Y\mathfrak{G}_v(X,Y)\big|_{X=Y}=\begin{cases}
1, & u=v,\\
0, & u\neq v.
\end{cases}$$
\end{Coro}
\textsc{Proof.} This is a standard computation
$$\mathfrak{G}_v(X,X)=\mathfrak{G}_v(X,Y)\bigg|_{Y=X}=\sum_{u\in W} \mathfrak{G}_u(\bullet,Y)\Pi_{u^{-1}}^\bullet \mathfrak{G}_v(X,\bullet)\bigg|_{\bullet=Y}\bigg|_{Y=X}.$$
Change the variables, we get the assertion. 
\qquad Q.E.D. \bigbreak

\term We can get more interesting combinatoricial identities. 

\begin{Th}[Star-Cauchy Formula=Theorem \ref{CauchyformulaforGrothendieck}]\label{StarCauchyFor}We have the following identity
$$\mathfrak{G}_w(X,Z)=\sum_{u*v=w}(-\beta)^{\ell(u)+\ell(v)-\ell(w)} \mathfrak{G}_v(X,Y)\mathfrak{G}_u(Y,Z)$$
\end{Th}
\textsc{Proof.} Apply the inversion formula \ref{inversionFo1} twice, 
$$\begin{array}{rl}
f(Z)&\displaystyle
=\sum_{w\in W} \mathfrak{G}_w(X,Z)\Pi_{w^{-1}}^X f(X)\\
& \displaystyle=\sum_{u\in W} \mathfrak{G}_u(Y,Z)\Pi_{u^{-1}}^Y f(Y)\\
& \displaystyle
=\sum_{u\in W} \mathfrak{G}_u(Y,Z)
\sum_{v\in W} \mathfrak{G}_v(X,Y) \Pi_{v^{-1}}^X\Pi_{u^{-1}}^X f(X)
\end{array}$$
Since $f$ is arbitrary, we can compare the coefficients of $\Pi_{w^{-1}}^Xf(X)$. 
\qquad Q.E.D. \bigbreak

\begin{Th}[Reduced Cauchy Formula]\label{RedCauchy}
In $R_T^\beta \otimes_{R_G^\beta}R_T^\beta \otimes_{R_G^\beta}R_T^\beta $, 
$$\begin{array}{l}
\displaystyle
\quad X^{\beta \rho}\mathfrak{G}_{w_0}(X,Y)=(-1)^{\ell(w_0)}Y^{\beta\rho}\mathfrak{G}_{w_0}(Y,X)\\[2ex]
\displaystyle
=\sum_{u\odot v=w_0}(-1)^{\ell(v)}Z^{\beta\rho}\mathfrak{G}_{v^{-1}}(Z,X)\mathfrak{G}_{u}(Z,Y). 
\end{array}$$
\end{Th}
\textsc{Proof.} 
We apply Corollary \ref{SymmetricTheorem} to
$f=\mathfrak{G}_{w_0}(Z,X), g=\mathfrak{G}_{w_0}(Z,Y)$. 
We get
$$\begin{array}{l}
\quad \displaystyle
\sum_{v\in W} (-1)^{\ell(v)} Z^{\beta \rho}\mathfrak{G}_{v^{-1}}(Z,X)\mathfrak{G}_{w_0v^{-1}}(Z,Y)\\
\displaystyle
=\sum_{w_0=u\odot v} (-1)^{\ell(v)} Z^{\beta \rho}\mathfrak{G}_{v^{-1}}(Z,X)\mathfrak{G}_{u}(Z,Y)
\end{array}$$
is symmetric in $Z$. Thus we can exchange $Z$ to $X$ or $Y$
$$\begin{array}{l}
\displaystyle
\quad \sum_{w_0=u\odot v} (-1)^{\ell(v)} Z^{\beta \rho}\mathfrak{G}_{v^{-1}}(Z,X)\mathfrak{G}_{u}(Z,Y)\\
\displaystyle=
\quad \sum_{w_0=u\odot v} (-1)^{\ell(v)} X^{\beta \rho}\mathfrak{G}_{v^{-1}}(X,X)\mathfrak{G}_{u}(X,Y)=X^{\beta \rho}\mathfrak{G}_{w_0}(X,Y)\\
\displaystyle=
\quad \sum_{w_0=u\odot v} (-1)^{\ell(v)} Y^{\beta \rho}\mathfrak{G}_{v^{-1}}(Y,X)\mathfrak{G}_{u}(Y,Y)=(-1)^{w_0}Y^{\beta \rho}\mathfrak{G}_{w_0}(Y,X)\\
\end{array}$$
This is what claimed in the proposition. 
\qquad Q.E.D. \bigbreak

\begin{Prop}\label{Mobiudinversionancient}We have 
$$X^{\beta\rho}\mathfrak{G}_{w}(X,Y)=\sum_{v\geq w}
(-1)^{\ell(v)} \beta^{\ell(v)-\ell(w)}Y^{\beta \rho}\mathfrak{G}_{v^{-1}}(Y,X).$$
\end{Prop}
\textsc{Proof.} 
We apply Corollary \ref{SymmetricTheorem} to
$f=\mathfrak{G}_{w_0}(Z,X), g=\mathfrak{G}_{w}^{\beta}(Z,Y)$. 
The element 
$$\begin{array}{l}
\quad \displaystyle
\sum_{v\in W} (-1)^{\ell(v)} Z^{\beta \rho}\mathfrak{G}_{v^{-1}}(Z,X)\pi_v^Z\mathfrak{G}_{w}(Z,Y)\\
\end{array}$$
is symmetric in $Z$. Thus we can exchange $Z$ to $X$ or $Y$
$$\begin{array}{l}
\quad \displaystyle
\sum_{v\in W} (-1)^{\ell(v)} Z^{\beta \rho}\mathfrak{G}_{v^{-1}}(Z,X)\pi_v^Z\mathfrak{G}_{w}(Z,Y)\\
\displaystyle=
\sum_{v\in W} (-1)^{\ell(v)} X^{\beta \rho}\mathfrak{G}_{v^{-1}}(X,X)\pi_v^X\mathfrak{G}_{w}(X,Y)=X^{\beta\rho}\mathfrak{G}_{w}(X,Y)\\
\displaystyle=
\sum_{v\in W} (-1)^{\ell(v)} Y^{\beta \rho}\mathfrak{G}_{v^{-1}}(Y,X)\pi_v^Z\mathfrak{G}_{w}(Z,Y)\bigg|_{Z=Y}\\
\displaystyle=\sum_{v\geq w}
(-1)^{\ell(v)} \beta^{\ell(v)-\ell(w)}Y^{\beta \rho}\mathfrak{G}_{v^{-1}}(Y,X)
\end{array}$$
Here we use the fact $\pi_u\mathfrak{G}_w$ is a constant if and only if $u\geq w$. 
\qquad Q.E.D. \bigbreak

\section{Generalized Billey formula}

\term Taking in $Y=uX$ in Theorem \ref{inversionFo1}, we get 

\begin{Th}\label{Localizationvalue}The value of the localization $\mathfrak{G}_w(X,uX)$ is determined by the following properties in $\NH^\beta(W)$
$$u=\sum_{w\in W} \mathfrak{G}_w(X,uX)\Pi_{w^{-1}},$$
for any $u\in W$.
\end{Th}

\def\TT{\mathbf{T}}
\def\L#1{\mathbf{L}_{#1}}
\def\ww{\underline{w}}
\def\hh{\mathbf{h}}

\term By the description of the localization map, Theorem \ref{GeoBilleyformula} is equivalent to 

\begin{Th}[Generalized Billey Formula]\label{GenBilleyFormula}
For the localizations of the generalized double $\beta$-Grothendieck polynomials, 
$$\mathfrak{G}_u(wx,x)=\sum_{w(J)=u}(-\beta)^{|J|-\ell(u)}\prod_{j\in J} d_j. $$
\end{Th}
The notation is introduced before Theorem \ref{GeoBilleyformula}. 

\term For example, we consider the case of $A$-type. 
Now, $a_i=\frac{x_i-x_{i+1}}{1-\beta x_{i+1}}=\beta^{-1}\big(1-\mbox{$\frac{X_{i+1}}{X_i}$}\big)$. 
Thus $\bar{a}_i=\frac{x_{i+1}-x_{i}}{1-\beta x_{i}}=\beta^{-1}\big(1-\mbox{$\frac{X_{i}}{X_{i+1}}$}\big)$.

Consider the case $w=s_2s_1$ and $u=s_2$, 
$$\begin{array}{c|c|c|c}\hline
w = s_2s_1& |J| & d_j& \Pi d_j\\\hline
s_2 = s_2\phantom{s_1} & 1 & \bar{a}_2\phantom{s_1}
& \beta^{-1}(1-\mbox{$\frac{X_2}{X_3}$})\\\hline
\end{array}$$
So $\mathfrak{S}_u(wx,x)=\beta^{-1}(1-\mbox{$\frac{X_2}{X_3}$})$. 
This coincides with our computation before. 

Similarly, for the case $w=s_1s_2$ and $u=s_2$, 
$$\begin{array}{c|c|c|c}\hline
w = s_1s_2& |J| & d_j& \Pi d_j\\\hline
s_2 = \phantom{s_1}s_2 & 1 & s_1\bar{a}_2
& \beta^{-1}(1-\mbox{$\frac{X_1}{X_3}$})\\\hline
\end{array}$$
So $\mathfrak{S}_u(wx,x)=\beta^{-1}(1-\mbox{$\frac{X_1}{X_3}$})$.

\term 
Consider $w=s_2s_1s_2$, and $u=s_2$. 
$$\begin{array}{c|c|c|c}\hline
w = s_2s_1s_2& |J| & d_j& \Pi d_j\\\hline
s_2 = s_2\phantom{s_1s_2} & 1 & \bar{a}_2\phantom{s_1s_2}
& \beta^{-1}(1-\mbox{$\frac{X_2}{X_3}$})\\\hline
s_2 = \phantom{s_2s_1}s_2 & 1 & s_2s_1\bar{a}_2
& \beta^{-1}(1-\mbox{$\frac{X_1}{X_2}$})\\\hline
s_2 = s_2\,*\,s_2 & 2 & \begin{array}{@{}l@{}}
\bar{a}_2\phantom{s_1s_2}\\s_2s_1\bar{a}_2\end{array}
& \beta^{-2}(1-\mbox{$\frac{X_2}{X_3}$})(1-\mbox{$\frac{X_1}{X_2}$})\\\hline
\end{array}$$
Thus 
$$\mathfrak{G}_u(wx,x)=\beta^{-1}(1-\mbox{$\frac{X_2}{X_3}$})+
\beta^{-1}(1-\mbox{$\frac{X_1}{X_2}$})
-\beta^{-1}(1-\mbox{$\frac{X_2}{X_3}$})(1-\mbox{$\frac{X_1}{X_2}$}).$$
This coincides with our computation before. 
On the other hand, we can take $w=s_1s_2s_1$, 
Then
$$\begin{array}{c|c|c|c}\hline
w = s_1s_2s_1& |J| & d_j& \Pi d_j\\\hline
s_2 = \phantom{s_1}s_2\phantom{s_1} & 1 & s_1\bar{a}_2\phantom{s_1}
& \beta^{-1}(1-\mbox{$\frac{X_1}{X_3}$})\\\hline
\end{array}$$
So
$$\mathfrak{G}_u(wx,x)=\beta^{-1}(1-\mbox{$\frac{X_1}{X_3}$}).$$
One can check that they are equal.

\term For $A$-type, above process can be interpolated to the language of pipe dreams when the decomposition of $w$ is reduced. 
For two pipe dreams $\pi_1$ and $\pi_2$, we write $\pi_1\leq \pi_2$ if they are the same at all the positions  where $\pi_2$ is not tiled by $\pipe$. 

Assume further that $\pi_2$ does not contain $\Bpipe$. Then we define the following weights for $(\pi_1,\pi_2)$.
$$\wt(\,\bump,\,\bump\,)=\wt(\,\bump,\,\pipe\,)=1.\qquad 
\begin{cases}
\wt(\,\pipe,\,\pipe\,)=\dfrac{x_{j}-x_i}{1-qx_i}=\beta^{-1}\big(1-\mbox{$\frac{X_j}{X_i}$}\big)\\[1ex]
\wt(\,\Bpipe,\,\pipe\,)=-\beta\frac{x_{j}-x_i}{1-qx_i}=-\big(1-\mbox{$\frac{X_j}{X_i}$}\big)
\end{cases}$$
where the two pipes in the second entry go to the $i$-th and $j$-th endings in $\pi_2$ with $i>j$ (NOT in $\pi_1$! ). 
Then define $\wt(\pi_1,\pi_2)$ by the product of weight at each position. 
Then for a fix one pipe dream $\pi_0$ consisting only $\pipe$ and $\bump$ with $w(\pi_0)=w$, and permutation $u$, 
$$\mathfrak{G}_u(wx,x)=\sum\nolimits_{\begin{subarray}{c}
{\pi\leq \pi_0}\\{w(\pi)=u}\end{subarray}}\wt(\pi,\pi_0). $$

\term For example, let $w=s_2s_3s_2$, and $u=s_3$. 
$$
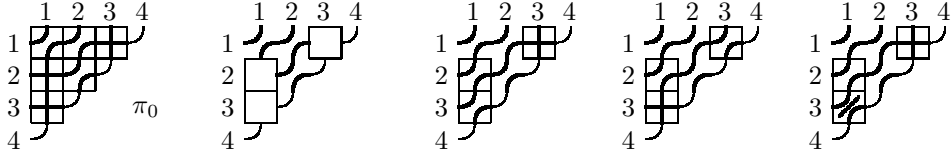
\begin{array}{l}
\labe{}\labe{1}\labe{2}\labe{3}\labe{4}\\
\labe{1}\bump\bump\pipe\abou\\
\labe{2}\pipe\bump\abou\\
\labe{3}\pipe\abou\hfill \pi_0\\
\labe{4}\abou
\end{array}\quad 
\begin{array}{l}
\labe{}\labe{1}\labe{2}\labe{3}\labe{4}\\
\labe{1}\Bbump\Bbump\empt\abou\\
\labe{2}\empt\Bbump\abou\\
\labe{3}\empt\abou\\
\labe{4}\abou
\end{array}\quad 
\begin{array}{l}
\labe{}\labe{1}\labe{2}\labe{3}\labe{4}\\
\labe{1}\Bbump\Bbump\pipe\abou\\
\labe{2}\bump\Bbump\abou\\
\labe{3}\bump\abou\\
\labe{4}\abou
\end{array}\begin{array}{l}
\labe{}\labe{1}\labe{2}\labe{3}\labe{4}\\
\labe{1}\Bbump\Bbump\bump\abou\\
\labe{2}\bump\Bbump\abou\\
\labe{3}\pipe\abou\\
\labe{4}\abou
\end{array}\begin{array}{l}
\labe{}\labe{1}\labe{2}\labe{3}\labe{4}\\
\labe{1}\Bbump\Bbump\pipe\abou\\
\labe{2}\bump\Bbump\abou\\
\labe{3}\Bpipe\abou\\
\labe{4}\abou
\end{array}$$
Thus 
$$\mathfrak{G}_u(wx,x)=\beta^{-1}(1-\mbox{$\frac{X_3}{X_4}$})+
\beta^{-1}(1-\mbox{$\frac{X_2}{X_3}$})
-\beta^{-1}(1-\mbox{$\frac{X_3}{X_4}$})(1-\mbox{$\frac{X_2}{X_3}$}).$$
On the other hand, we can pick
$$
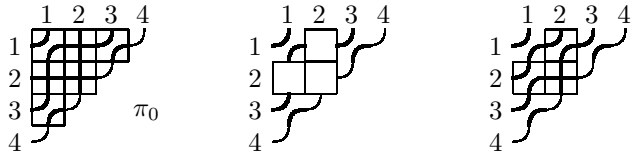
\begin{array}{l}
\labe{}\labe{1}\labe{2}\labe{3}\labe{4}\\
\labe{1}\bump\pipe\bump\abou\\
\labe{2}\pipe\pipe\abou\\
\labe{3}\bump\abou\hfill \pi_0\\
\labe{4}\abou
\end{array}\qquad 
\begin{array}{l}
\labe{}\labe{1}\labe{2}\labe{3}\labe{4}\\
\labe{1}\Bbump\empt\Bbump\abou\\
\labe{2}\empt\empt\abou\\
\labe{3}\Bbump\abou\\
\labe{4}\abou
\end{array}\qquad 
\begin{array}{l}
\labe{}\labe{1}\labe{2}\labe{3}\labe{4}\\
\labe{1}\Bbump\bump\Bbump\abou\\
\labe{2}\bump\pipe\abou\\
\labe{3}\Bbump\abou\\
\labe{4}\abou
\end{array}$$
Thus 
$$\mathfrak{G}_u(wx,x)=\beta^{-1}(1-\mbox{$\frac{X_2}{X_4}$}). $$
This coincides with the Buch--Rim\'anyi formula introduced in \cite{buch2004specializations}.
See \cite{fan2020note} for the combinatorial model of it. 

\term[Proof of Theorem \ref{GenBilleyFormula}]
For any $f(x)\in R_T^\beta$, we denote $\L{f(x)}:\NH^\beta(W)\to \NH^\beta(W)$ the left multiplication by $f(x)$. 
Denote $\TT_w:\NH^\beta(W)\to \NH^\beta(W)$ the right multiplication by $\Pi_{w^{-1}}^X$, and $\TT_i=\TT_{s_i}$ for simplicity. 
We also denote $a_i=\frac{1-e^{-\beta\alpha_i}}{\beta}$. 
We denote $\hh_i(x)=1+\L{x}\TT_i$. Then 
$$\begin{array}{rl}
ws_i& = w\big(1+\frac{1-e^{-\beta\alpha_i}}{\beta}\Pi_i\big)=w(1+a_i\Pi_i)\\[1ex]
& =w+(wa_i)w\Pi_i
=(1+\L{wa_i}T_i)(w)=\hh_i(wa_i)(w). 
\end{array}$$
Since $\L{x}$ commutes with $\TT_i$'s, $\hh_i$'s satisfy the Yang--Baxter equation in $A$-type above (while we will not use this fact). 
Recall that $w=s_{i_1}s_{i_2}\cdots s_{i_r}$. 
Let us denote $\ww^{(j)}=s_{i_1}\cdots s_{i_j}$. 
Then 
$$\begin{array}{rl}
w & = \ww^{(r-1)} s_{i_r}=\hh_{i_r}(w^{(r-1)}a_{i_r})(\ww^{(r-1)})\\[1ex]
& =\big(\hh_{i_r}(\ww^{(r-1)}a_{i_r})\circ \hh_{i_{r-1}}(\ww^{(r-2)}a_{i_{r-1}})\big) (\ww^{(r-2)})\\[1ex]
& =\big(\hh_{i_r}(\ww^{(r-1)}a_{i_r})\circ \cdots \circ \hh_{i_1}(\ww^{(0)}a_{i_1})\big)(\operatorname{id})
\end{array}$$
Assume 
$$\hh_{i_r}(\ww^{(r-1)}a_{i_r})\circ \cdots \circ \hh_{i_1}(\ww^{(0)}a_{i_1})=\sum_{u\in W} \L{c^u_w(x)}\TT_u$$
Then $w=\sum c^u_w(x) \Pi_{u^{-1}}$. 
By Theorem \ref{Localizationvalue}, 
$c^u_w(x)=\mathfrak{G}_u(x,wx)$. 

Under the notation of the previous subsection, we get the following identity
$$h_{i_r}(\ww^{(r-1)}a_{i_r})\cdots  h_{i_1}(\ww^{(0)}a_{i_1})
=\sum_{u\in W}\mathfrak{G}_u(x,wx)T_u, $$
By applying $w^{-1}$ both side and replace $w$ by $w^{-1}$, we get
$$h_{i_1}(\ww^{(0)}\bar{a}_{i_1})\cdots  h_{i_r}(\ww^{(r-1)}\bar{a}_{i_r})
=\sum_{u\in W} \mathfrak{G}_u(wx,x)T_u. $$
Here $\bar{a}_i=\ominus_{\beta} a_i:=\frac{-a_i}{1-\beta a_i}$. 
That is, 
$$h_{i_1}(d_1)\cdots  h_{i_r}(d_r)
=\sum_{u\in W} \mathfrak{G}_u(wx,x)T_u. $$
The expansion of the left-hand-side gives the formula in Theorem \ref{GenBilleyFormula}. 
\qquad Q.E.D. \bigbreak


\section{Dual Grothendieck Polynomials}

\term We define the \emph{dual Grothendieck polynomial}
$$\mathfrak{g}_{w^{-1}}(X,Y)=\Pi_{w^{-1}w_0}^Y\mathfrak{G}_{w_0}(X,Y). $$
Note that $\lim_{\beta \to 0}\mathfrak{g}_{w}(X,Y)=\mathfrak{S}_{w}(x,y)$. 

\term The following identity is the analogy of the second identity of Theorem \ref{doubleidentities}. 

\begin{Th}We have  
$$X^{\beta\rho} \mathfrak{G}_{w}(X,Y)=(-1)^{\ell(w)} Y^{\beta\rho}\mathfrak{g}_{w^{-1}}(Y,X).$$
In particular, 
$$u\leq w\Longrightarrow\mathfrak{g}_{w}(uX,X)=0,\qquad \mathfrak{g}_{\operatorname{id}}(X,Y)=X^{\beta\rho}Y^{-\beta\rho}.$$
\end{Th}
\textsc{Proof.} 
We can compute
$$\begin{array}{rl}
X^{\beta\rho} \mathfrak{G}_{w}(X,Y)&
= X^{\beta\rho} \pi^X_{w^{-1}w_0}\mathfrak{G}_{w_0}(X,Y)\\
&= (-1)^{\ell(w_0)}X^{\beta\rho} \pi^X_{w^{-1}w_0}X^{-\beta\rho}Y^{\beta\rho}\mathfrak{G}_{w_0}(Y,X)\\
&= (-1)^{\ell(w_0)}X^{\beta\rho} \pi^X_{w^{-1}w_0}D\big(\overline{Y^{\beta\rho}\mathfrak{G}_{w_0}(Y,X)}\big)\\
&= (-1)^{\ell(w)}X^{\beta\rho} D\big(\ppi^X_{w^{-1}w_0} \overline{Y^{\beta\rho}\mathfrak{G}_{w_0}(Y,X)}\big)\\
&= (-1)^{\ell(w)}\Pi^X_{w^{-1}w_0}Y^{\beta\rho}\mathfrak{G}_{w_0}(Y,X)
=(-1)^{\ell(w)} Y^{\beta\rho}\mathfrak{g}_{w^{-1}}(Y,X).
\end{array}$$
The localization condition is trivial. 
\qquad Q.E.D. \bigbreak

\term 
We can rewrite the reduced Cauchy formula \ref{RedCauchy} and Proposition \ref{Mobiudinversionancient}. 

\begin{Coro}[Reduced Cauchy Formula]\label{RedCauchy2}
In $K$,  
$$\mathfrak{G}_{w_0}(X,Y)=\sum\nolimits_{u\odot v=w_0}\mathfrak{g}_{v}(X,Z)\mathfrak{G}_{u}(Z,Y). $$
\end{Coro}

\begin{Coro}[M\"obius Inversion]\label{MobiusInversion}
We have 
$$\begin{array}{c}
\displaystyle
\mathfrak{G}_{w}(X,Y)=\sum_{v\geq w}
\beta^{\ell(v)-\ell(w)}\mathfrak{g}_{v}(X,Y),\\[2ex]
\displaystyle
\mathfrak{g}_{w}(X,Y)=\sum_{v\geq w}
(-\beta)^{\ell(v)-\ell(w)}\mathfrak{G}_{v}(X,Y).
\end{array}$$
\end{Coro}

\term Dual to Theorem \ref{inversionFo1} and Corollary \ref{DualFor1} above, we have 

\begin{Th}[Inversion Formula]\label{inversionFo2}For any $f\in R_T$, we have the following identity in $K$
$$f(X)=\sum_{w\in W} \mathfrak{g}_w(X,Y)\pi_w^Y f(Y)$$
\end{Th}
\textsc{Proof.} Note that the right-hand-side is 
$$\Pi_{ww_0}^Y\mathfrak{G}_{w_0}(X,Y)\pi_w^Yf(Y)$$
which is symmetric in $Y$ by Corollary \ref{SymmetricTheorem2}, so that we can take $Y=X$. 
\qquad Q.E.D. \bigbreak

\begin{Coro}[Dual Basis]\label{DualFor2}The operator $\pi_{w}$ and $\mathfrak{g}_w$ are dual to each other, that is, 
$$\pi_{u}^X\mathfrak{g}_v^{\beta}(X,Y)\big|_{Y=X}=\begin{cases}
1, & u=v,\\
0, & u\neq v.
\end{cases}$$
\end{Coro}

\begin{Th}[Star-Cauchy Formula]\label{StarCauchyFor2}
We have the following identity
$$\mathfrak{g}_w(X,Z)=\sum\nolimits_{u*v=w}\beta^{\ell(u)+\ell(v)-\ell(w)} \mathfrak{g}_v(X,Y)\mathfrak{g}_u(Y,Z)$$
\end{Th}
\textsc{Proof.} The proof is the same to Theorem \ref{StarCauchyFor}. 
\qquad Q.E.D. \bigbreak

\section{Appendix: Pipe Dream}

\term We will shortly review the combinatorics in $A$-type which provides examples for the next subsection.
In the case $G=\GL_n$. We recognize $\Lambda=\mathbb{Z}t_1\oplus \cdots \oplus \mathbb{Z}t_n$. 
In cohomology, we denote $x_i=t_i$ (not confusing with the K-theory notation). 
The Demazure operator 
$$\partial_if=\frac{f-s_if}{x_i-x_j}. $$
Other than the Euler class, we have another choice \cite{lascoux2007anneau} 
$\mathfrak{S}_{w_0}(x)=x_1^{n-1}x_2^{n-2}\cdots x_{n-1}$, and 
$\mathfrak{S}_{w_0}(x,y)=\prod_{i+j\leq n} (x_i-y_j)$. 
This choice is stable. 

\term 
We denote $X_i=e^{\beta t_i}\otimes 1$, $Y_i=1\otimes e^{\beta t_i}$, and $x_i=\frac{1-e^{-\beta t_i}}{\beta}\otimes 1$, 
$y_i=1\otimes \frac{1-e^{-\beta t_i}}{\beta}$. 
Then in terms of $X_i$, 
$$\pi_i f= \beta\frac{f-e^{-\beta\alpha_i}s_if}{1-e^{-\beta\alpha_i}}
=\beta\frac{f-\frac{X_{i+1}}{X_i}s_if}{1-\frac{X_{i+1}}{X_i}}
=\beta\frac{X_if-X_{i+1}s_if}{X_i-X_{i+1}}
=\beta\cdot  \partial_i^X(X_if).$$
In terms of $x_i$, 
$$\pi_i f=\beta\frac{\frac{1}{1-\beta x_i}f-\frac{1}{1-\beta x_{i+1}}s_if}{\frac{1}{1-\beta x_i}-\frac{1}{1-\beta x_{i+1}}}
=\frac{(1-\beta x_{i+1})f-(1-\beta x_{i})s_if}{x_i-x_{i+1}}
=\partial_i((1-\beta x_{i+1})f)
.$$
We have the stable choice \cite{fomin1993yang} (after replacing $\beta$ by $-\beta$) 
$$\mathfrak{G}_{w_0}(X,Y)=\beta^{-\ell(w_0)}\prod_{i+j\leq n} \big(1-\mbox{$\frac{Y_i}{X_i}$}\big)= 
\prod_{i+j\leq n} \frac{x_i-y_j}{1-\beta y_j}. $$

\term 
We denote $h_i(x)=(1+xT_i)$, where $T_i$'s satisfy 
$T_i^2=-\beta T_i$ and the Braid relations. 

In $A$-type, one can check the following \emph{Yang--Baxter equations}
$$\begin{array}{c}
h_i(x)h_i(y)=h_i(x\oplus_\beta y)\\
h_i(x)h_j(y)=h_j(y)h_i(x)\qquad  |i-j|\geq 2\\
h_i(x)h_{i+1}(x\oplus_\beta y) h_i(y)=
h_{i+1}(y)h_i(x\oplus_\beta y)h_{i+1}(x)
\end{array}$$
Here the variables $x$ and $y$ commute with $T_i$'s. 
There are also Yang--Baxter equations for non-simply-laced cases, see \cite{kirillov2015double} and \cite{billey1999kostant}. 

\term 
We consider the generating function 
$$\mathfrak{G}^\beta(x,y)=\sum \mathfrak{G}_{w}^\beta(x,y)T_{w}^\beta. $$
It is amazing that it factors into 
$$\begin{array}{ccccc}
h_{n-1}(x_1\ominus y_{n-1})&h_{n-2}(x_1\ominus y_{n-2})&\cdots& h_1(x_1\ominus y_2)&h_1(x_1\ominus y_1)\\
&h_{n-1}(x_2\ominus y_{n-2})&\cdots &h_3(x_2\ominus y_2)&h_2(x_1\ominus y_1)\\
&&\ddots &\vdots &\vdots \\
&&& h_{n-1}(x_{n-2}\ominus y_2) & h_{n-2}(x_{n-2}\ominus y_1)\\
&&&& h_{n-1}(x_{n-1}\ominus y_1)
\end{array}$$
where $x\ominus_{\beta} y:= \frac{x-y}{1-\beta y}$, see \cite{fomin1996yang}, it is the proof for $\beta=0$, but the proof works in general. 
This observation originally appeared in \cite{fomin1994schubert}, the case $\beta=0, y_i=0$. 
One can prove (see \cite{fomin1996yang})
$$\mathfrak{G}^\beta(y,z)\mathfrak{G}^\beta(x,y)=\mathfrak{G}^\beta(x,z)$$
which implies the star-Cauchy formula \ref{StarCauchyFor} for $A$- type.

\term The expansion of this expression gives the model known as \emph{pipe dreams} \cite{bergeron1993rc} (the case $\beta=0$ and $y=0$), and \cite{knutson2019schubert} (the case $\beta=0$).  
To be exact, a pipe dream is a tiling of the type $(n-1,\ldots,1,0)$ tableau by three kinds of tiles
$\pipe$, $\bump$ and $\Bpipe$ in 
$\begin{array}{l}
\labe{}\labe{1}\labe{2}\labe{\cdots}\labe{\cdots}\labe{n}\\
\labe{1}\empt\empt\empt\empt\abou\\
\labe{2}\empt\empt\empt\abou\\
\labe{\vdots}\empt\empt\abou\\
\labe{\vdots}\empt\abou\\
\labe{n}\abou\\
\end{array}$ such that 
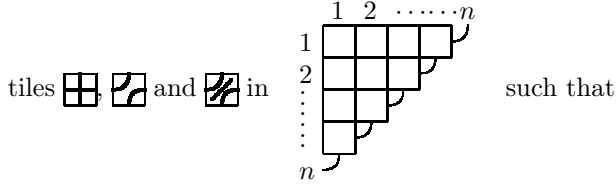
\begin{itemize}
    \item any pair of pipes intersect at most once and
    \item the two pipes in any appearance $\Bpipe$ intersects at the northeast of it.
\end{itemize}
For each pipe dream $\pi$, one can read a permutation $w(\pi)$ from the left boundary to the up boundary. 
We define the weight of tiles at $(i,j)$-position by 
$$\wt(\,\bump\,)=1,\qquad
\begin{cases}
\wt(\,\pipe\,)=\dfrac{x_i-y_j}{1-\beta y_j}=\beta^{-1}\big(1-\mbox{$\frac{Y_i}{X_i}$}\big),\\[1ex]
\wt(\,\Bpipe\,)=-\beta \dfrac{x_i-y_j}{1-\beta y_j}=-\big(1-\mbox{$\frac{Y_i}{X_i}$}\big).
\end{cases} $$
The weight $\wt(\pi)$ of a pipe dream $\pi$ is defined to be the product of weights of all tiles. 
Then 
$$\mathfrak{G}_w(x,y)=\sum_{w(\pi)=w} \wt(\pi). $$

\term 
For example, for the permutation $w=\binom{123}{132}\in\mathfrak{S}_3$, 
$$\begin{array}{ccccc}
&
\begin{array}{l}
\labe{}\labe{1}\labe{2}\labe{3}\\
\labe{1}\bump\bump\abou\\
\labe{2}\pipe\abou\\
\labe{3}\abou
\end{array}
& 
\begin{array}{l}
\labe{}\labe{1}\labe{2}\labe{3}\\
\labe{1}\bump\pipe\abou\\
\labe{2}\bump\abou\\
\labe{3}\abou
\end{array}
&
\begin{array}{l}
\labe{}\labe{1}\labe{2}\labe{3}\\
\labe{1}\bump\pipe\abou\\
\labe{2}\Bpipe\abou\\
\labe{3}\abou
\end{array}
&
\begin{array}{l}
\labe{}\labe{1}\labe{2}\labe{3}\\
\labe{1}\bump\Bpipe\abou\\
\labe{2}\pipe\abou\\
\labe{3}\abou
\end{array}
\\\\
\wt&\frac{x_2-y_1}{1-\beta y_1}& 
\frac{x_1-y_2}{1-\beta y_2} & 
-\beta \big(\frac{x_1-y_2}{1-\beta y_2}\big)
\big(\frac{x_2-y_1}{1-\beta y_1}\big)& \text{Not allowed}\\[2ex]
\wt&\beta^{-1}\big(1-\mbox{$\frac{Y_1}{X_2}$}\big)&
\beta^{-1}\big(1-\mbox{$\frac{Y_2}{X_1}$}\big)&
-\beta^{-1}\big(1-\mbox{$\frac{Y_2}{X_1}$}\big)\big(1-\mbox{$\frac{Y_1}{X_2}$}\big)& \text{Not allowed}
\end{array}$$
Thus 
$$\begin{array}{rl}
\mathfrak{G}_w(x,y)& =\frac{x_2-y_1}{1-\beta y_1}+
\frac{x_1-y_2}{1-\beta y_2}
-\beta \big(\frac{x_1-y_2}{1-\beta y_2}\big)
\big(\frac{x_2-y_1}{1-\beta y_1}\big)\\[2ex]
&=\beta^{-1}\big(1-\frac{Y_1}{X_2}\big)+\beta^{-1}\big(1-\frac{Y_2}{X_1}\big)-\beta^{-1}
\big(1-\frac{Y_2}{X_1}\big)\big(1-\frac{Y_1}{X_2}\big).
\end{array}$$

\bibliography{references}

\end{document}